\documentclass[11pt,a4wide]{article}
\usepackage[latin,english]{babel}
\usepackage[utf8]{inputenc}
\usepackage[T1]{fontenc}
\usepackage{eufrak}
\usepackage{geometry}
\usepackage{amsmath}
\usepackage{enumitem}
\usepackage{amsfonts}
\usepackage{amssymb}
\usepackage{amsbsy}
\usepackage{amsthm}
\usepackage{makeidx}
\usepackage{mathtools}
\usepackage{mathabx, mathrsfs, dsfont}
\usepackage{cases}
\usepackage{braket}
\usepackage[toc,page]{appendix}
\usepackage{dirtytalk}
\usepackage{pdfsync}
\usepackage{graphicx}
\graphicspath{ {./Figures/} } 
\usepackage{epstopdf}
\usepackage{todonotes}
\usepackage{fancyhdr}
\usepackage{math}
\usepackage{cite}
\usepackage{color}
\usepackage{tikz}
\usepackage{subcaption}
\usepackage{bbm}
\usepackage{authblk}
\usepackage[colorlinks=true, pdfstartview=FitV, urlcolor=blue, citecolor=red, linkcolor=blue,pdfencoding=unicode,unicode=true,psdextra]{hyperref}
\usepackage{marginnote}

\mathtoolsset{showonlyrefs}

\usetikzlibrary{matrix,decorations.markings}
\usetikzlibrary{decorations.markings}
\usetikzlibrary{quotes,angles}

\usepackage{esint}

\numberwithin{equation}{section}

\title{A Short Survey of the Well-posedness of the Two-dimensional Burgers' Equation}
\author{Xiang Zhang}
\author{Yule Sun}
\author{Shuhan Xie}
\affil{School of Mathematical Science, Tongji University}

\date{\today}

\begin{document}
\maketitle

\abstract{
In this paper, we establish the existence and uniqueness of solutions to the two-dimensional Burgers equation using the framework of infinite-dimensional dynamical systems. The two-dimensional Burgers equation, which models the interplay between nonlinear advection and viscous dissipation, is given by:
$$
u_{t} + u \cdot \nabla u = \nu \Delta u + f,
$$
where $ u = (u_1, u_2) $ is the velocity field, $ \nu > 0 $ is the viscosity coefficient, and $ f $ represents an external force.\cite{e2000} We primarily employed Galerkin method to transform the partial differential equation into an ordinary differential equation. In addition, by employing Sobolev spaces, energy estimates, and compactness arguments, we rigorously prove the existence of global solutions and their uniqueness under appropriate initial and boundary conditions.
}

\section{Introduction And Main Result}
\paragraph{}
We investigate the 2D viscous Burgers equation
\begin{equation}
    u_{t} + u \cdot \nabla u = \nu \Delta u + f,x \in \Omega
\end{equation}
where $ \Omega $ is a open bounded set in $ \mathbb{R}^{2} $ with boundary $ \Gamma $ and $ \Omega $ is smooth enough. The external force $ f $ in equation is independent of time $ t $. With the initial data
\begin{equation}
    u(0,x,y) = u_{0}(x,y)
\end{equation} 
The equation are supplemented with a boundary condition.Two cases will be considered:
\\\textit{The nonslip boundary condition.} The boundary $ \Gamma$ is solid and at rest;thus $$ u=0 \quad \text{on} \ \Gamma$$
\textit{The space-periodic case.} Here $ \Omega = (0,L_{1}) \times (0,L_{2}) $ and $ u, p $ and the first derivatives of $ u $ are $ \Omega\text{-}periodic. $
\paragraph{}
First, we consider a Hilbert space $ H = L^{2}(\Omega)$ and its closed subspace $ V = H_{0}^{1}(\Omega) $. We prove that $ V \subset H \subset V^{*} $ , $ V $ is dense in $ H $ and $ H $ is dense in $ V^{*} $. $ V^{*} $ is the dual space of $ V $. Then we denote by A the linear unbounded operator in $ H $:
\begin{equation}
    (Au,v)_{L^{2}(\Omega)} = (u,v)_{H_{0}^{1}(\Omega)}, \quad \forall u,v \in V
\end{equation}
It can be concluded that the operator A is an isomorphism from $ V $ onto $ V^{*} $ by analyzing the properties of the operator A. The domain of A in $ H $ is denoted by $ D(A) $, endowed with the norm $ \Vert Au \Vert_{H} $.
\paragraph{}
The weak form of the Burgers equation is obtained by multiplying by a test function $ \varphi $ in V and integrating over $ \Omega $:
\begin{equation}
    (P) \left\{
    \begin{aligned}
        \frac{d}{dt}u + vAu + B(u) &= f\\
        u(0,x,y) &= u_{0}(x,y)
    \end{aligned}
    \right.
\end{equation}
where B is the bilinear operator from $ V \times V $ into $ V^{*} $ defined by
\begin{flalign}
    (B(u,v),w) &= b(u,v,w), \quad \forall u,v,w \in V. \\
    where \quad b(u,v,w) &= \Sigma_{i,j=1}^{n} \int_{\Omega} u_{i} \cdot \partial_{x_{i}} v_{j} \cdot w_{j} dx.
\end{flalign}
Based on the spectral properties of the operator A, we can apply the Faedo-Galerkin method to transform the equation into a finite-dimensional case. In fact, in this scenario, the equation becomes an ordinary differential equation, which is easy to prove. We denote the solution of the equation under the m-dimensional case by $ u_{m} $. The original equation becomes
\begin{equation}
    \frac{d}{dt}u_{m} + vAu_{m} + P_{m}B(u_{m}) = P_{m}f
\end{equation}
where $ P_{m} $ is the projector in H on the m-dimensional space.
By using the energy-type equality, we obtain that $ u_{m} \in L^{\infty}(0,T;H) \cap L^{2}(0,T;V) \quad \forall T>0 $. Hence, $ Au_{m} \in L^{2}(0,T;V^{*}) $. Then we can prove that $ \partial_{t}u_{m} $ is bounded in $ L^{2}(0,T;V^{*}) $ by illustrating the boundedness of $ P_{m}B(u_{m}) $. In the end, the compactness theorem indicates there exists a function u in $ L^{2}(0,T;V^{*}) $ such that $ u_{m} $ strongly converges to $ u $. With two lemmas, we can prove $ u $ is the unique solution of our equation.
\begin{theorem}
    There exists a unique solution u of $ (P) $ satisfying
    \begin{equation}
        u \in C([0,T];H) \cap L^{2}(0,T;V), \quad \forall T>0.
    \end{equation}
    and the solution depends continuously on $ f $ and initial value $ u_{0} $.
    Furthermore, if $ u_{0} \in V $, then
    \begin{equation}
        u \in C([0,T];V) \cap L^{2}(0,T;D(A))
    \end{equation}
\end{theorem}

\section{Proof}
\subsection{The weak form and finite-dimensional case}
\paragraph{}
To prepare for the Faedo-Galerkin approximation, we need to prove the spectral properties of A to transform the equation into a finite dimensional case. Before the proof, we need to illustrate $ V \subset H \subset V^{*} $ and $ V $ is dense in $ H $, $ H $ is dense in $ V^{*} $. It is obvious that $ C_{0}^{\infty}(\Omega) $ is dense in $ L^{2}(\Omega) $ and $ H_{0}^{1}(\Omega) $. $ H_{0}^{1}(\Omega) \subset L^{2}(\Omega) $, which indicates $ V $ is dense in $ H $. $ H $ is a Hilbert space, so it is convenient to identify $ H $ to its dual $ H^{*} $ due to the Riesz representation theorem. Thus, we can identify $ H^{*} $ to a dense subspace of $ V^{*} $, so we obtain
$$
V \subset H \equiv H^{*} \subset V^{*}
$$
Then we embark on the proof of the spectral property. Firstly, we introduce a useful lemma:
\begin{lemma}[Lax-Milgram lemma]
If a is a bilinear continuous coercive form on V then A is an isomorphism from $ V $ onto $ V^{*} $. The coercive condition means there exists $ \alpha>0,a(u,u)\geq \alpha \Vert u \Vert_{V}^{2}, \quad \forall u \in V. $
\end{lemma}
Obviously, $ a(u,u) = \Vert u \Vert_{V} $, satisfies the condition when $ \alpha = 1 $, then $ A $ is an isomorphism from $ V $ onto $ V^{*} $. $ A $ is self-adjoint, because $ \langle Au,v \rangle = \langle Av,u \rangle = a(u,v), \quad \forall u,v \in V.$ The compact embedding theorem tell us that the injection of $ V $ in $ H $ is compact. In this case, the inverse $ A^{-1} $ can be considered as a self-adjoint compact operator in $ H $, hence, with the spectral theory:
\begin{theorem}[The spectral theorem]
    Suppose $ T $ is a compact symmetric operator on a Hilbert space $ H $. Then there exists an orthonormal basis $ \{ \varphi_{k} \}_{k=1}^{\infty} $ of $ H $ that consists of eigenvectors of $ T $. Moreover, if $ T\varphi_{k} = \lambda_{k}\varphi_{k} $, then $ \lambda_{k} \in \mathbb{R} $ and $ \lambda \rightarrow 0 $ as $ k \rightarrow \infty $. 
\end{theorem}
We infer that there exists a complete orthonormal family of H, $ \{w_{j}\}_{j\in \mathbb{N}} $ made of eigenvectors of $ A $:
$$
A^{-1}w_{j} = \mu_{j}w_{j}, \quad \forall j \in \mathbb{N},\ and\ \mu_{j} \rightarrow 0 \quad as \quad j \rightarrow \infty.
$$
It is clear that $ w_{j} \in D(A), \forall j \in \mathbb{N}$, and setting $ \lambda_{j} = \mu_{j}^{-1} $ we obtain:
$$
Aw_{j} = \lambda_{j}w_{j}, \quad \forall j \in \mathbb{N},\ and\ \lambda_{j} \rightarrow \infty \quad as \quad j \rightarrow \infty.
$$
Now, we have proven the spectral property and able to consider the problem projected on the m-dimensional space spanned by $ \{w_{j}\}_{j=1}^{m} $. In this case, the problem has been transformed into:
\begin{equation}
    \left\{
    \begin{aligned}
        \frac{d}{dt}u_{m} + vAu_{m} + P_{m}B(u_{m}) &= P_{m}f\\
        u_{m}(0,x,y) &= P_{m}u_{0}(x,y)
    \end{aligned}
    \right.
\end{equation}
where $ P_{m} $ is the projector.

\subsection{The Garlerkin approximation}
\paragraph{}
It can be readily observed that the original equation has been reformulated as an ordinary differential equation. By virtue of the Picard–Lindelöf existence and uniqueness theorem, there exists a constant $ T>0 $ such that the solution $ u_{m} $ exists uniquely on the interval $ [0,T) $.
\paragraph{}
In the subsequent analysis, we shall employ the energy method to proceed with the proof. The first energy-type inequality is obtained by taking the scalar product of the equation with $ u $, where $ u \in V $.
\begin{flalign}
    (\frac{d}{dt}u,u) + v(Au,u) + (B(u),u) &= (f,u)\\
    \Rightarrow \frac{1}{2} \frac{d}{dt} \Vert u \Vert_{L^{2}(\Omega)}^{2} + v\Vert u \Vert_{V}^{2} + (B(u),u) &\leq \Vert f \Vert_{L^{2}(\Omega)} \Vert u \Vert_{L^{2}(\Omega)} 
\end{flalign}
we see that $ (B(u),u) = \int_{\Omega} u \cdot Ju \cdot u\ dxdy = 0 $ by integration by parts. Furthermore, by invoking the Poincaré inequality, we deduce that $ \Vert u \Vert_{L^{2}(\Omega)} \leq \lambda_{1}^{-\frac{1}{2}} \Vert u \Vert_{V} $, $ \lambda_{1} $ is the first eigenvalue of $ A $. Substituting this into the above equation, we obtain that
\begin{flalign}
\frac{d}{dt} \Vert u \Vert_{L^{2}(\Omega)}^{2} + v\lambda_{1}\Vert u \Vert_{V}^{2} \leq \frac{1}{v\lambda_{1}} \Vert f \Vert_{L^{2}(\Omega)}^{2}
\end{flalign}
Using the classical Gronwall lemma, we obtain
\begin{flalign}
\Vert u(t) \Vert_{L^{2}(\Omega)}^{2} \leq \Vert u_{0} \Vert_{L^{2}(\Omega)}^{2} e^{-v\lambda_{1}t} + \frac{1}{v^{2}\lambda_{1}^{2}} \Vert f \Vert_{L^{2}(\Omega)}^{2}(1-e^{v\lambda_{1}t})
\end{flalign}
which demonstrates $ \Vert u(t) \Vert_{L^{2}(\Omega)}^{2} $ admits a uniform upper bound for all $ t $, hence $ u_{m}(t) $ exists uniquely on the interval $ [0,\infty) $.
Then we can easily get two conclusions
\begin{flalign}
u_{m} &\in L^{\infty}(0,T;H) \cap L^{2}(0,T;V) \quad \forall T>0 \\
Au_{m} &\in L^{2}(0,T;V^{*})
\end{flalign}

\paragraph{}
Now, we shift our focus to $ B(u_{m}) $:
\begin{flalign}
\forall \varphi \in V, \Vert \varphi \Vert_{V} &= 1 \\
\vert (B(u_{m},\varphi)) \vert \leq \int_{\Omega} \vert u_{m} \vert \vert \nabla u_{m} \vert \vert \varphi \vert dxdy &\leq \Vert \nabla u_{m} \Vert_{L^{2}(\Omega)} \Vert \varphi \Vert_{L^{2}(\Omega)} \Vert u_{m} \Vert_{L^{\infty}(\Omega)}   
\end{flalign}
Then, we introduce the Agmon's inequality\cite{agmon1965}:
\begin{lemma}[Agmon's inequality]
Assuming that $ \Omega \subset \mathbb{R}^{n} $ is of class $ C^{n} $, there exists a constant c depending only on $ \Omega $ such that:
\begin{equation}
    \Vert u \Vert_{L^{\infty}(\Omega)} \leq \left\{
    \begin{aligned}
        &c \Vert u \Vert_{H^{(n/2)-1}(\Omega)}^{1/2} \Vert u \Vert_{H^{(n/2)+1}(\Omega)}^{1/2}, \forall u \in H^{(n/2)+1}(\Omega) \quad if \ n \ is \ even, \\
        &c \Vert u \Vert_{H^{(n-1)/2}(\Omega)}^{1/2} \Vert u \Vert_{H^{(n+1)/2}(\Omega)}^{1/2}, \forall u \in H^{(n+1)/2}(\Omega) \quad if \ n \ is \ odd.
    \end{aligned}
    \right.
\end{equation}
\end{lemma}
By Agmon's inequality, we have:
\begin{flalign}
\Vert u_{m} \Vert_{L^{\infty}(\Omega)} \leq c\Vert u_{m} \Vert_{H^{0}(\Omega)}^{1/2} \Vert u_{m} \Vert_{H^{1}(\Omega)}^{1/2} = c\Vert u_{m} \Vert_{L^{2}(\Omega)}^{1/2} \Vert u_{m} \Vert_{H^{1}(\Omega)}^{1/2}
\end{flalign}
Hence, $ \vert (B(u_{m},\varphi)) \vert \leq c \Vert u_{m} \Vert_{V}^{3/2} \Vert u_{m} \Vert_{L^{2}(\Omega)}^{1/2} $, by the boundedness of $ u_{m} $, we obtain $ \Vert B(u_{m}) \Vert_{V} $ is bounded, thus $ B(u_{m}),P_{m}B(u_{m}) $ are bounded in $ L^{2}(0,T;V^{*}) $. Based on the equation we need to solve, we obtain the boundedness of $ \frac{du_{m}}{dt} $ in $ L^{2}(0,T;V^{*}) $. 

\paragraph{}
Next, we introduce two theorems firstly.
\begin{theorem}[Banach-Alaoglu Theorem]
Let $ X $ be a normed space, and let $ X^{*} $ be its dual space. Then, the closed unit ball $ B^{*} = \{ f \in X^{*}: \Vert f \Vert \leq 1 \} $
is compact in the weak-* topology.
\end{theorem}
\begin{theorem}[Compactness Theorem]
    If $ \Omega $ is a bounded domain in $ \mathbb{R}^{n} $ with smooth boundary, then the embedding of $ H_{0}^{1}(\Omega) $ into $ L^{2}(\Omega) $ is compact.
\end{theorem}
Due to the Banach-Alaoglu Theorem, there exists a sequence $ u_{m_{k}} $ that converges weakly to $ u $ in $ L^{2}(0,T;V) $ and weakly-star to u in $ L^{\infty}(0,T;H) $. Without loss of generality, we denote it as $ u_{m} $. Hence, $ u \in L^{\infty}(0,T;H) \cap L^{2}(0,T;V) $. Similarly, $ \frac{du_{m}}{dt} $ converges weakly in $ L^{2}(0,T;V^{*}) $. Due to the compactness theorem, $ V=H_{0}^{1}(\Omega) $ is compactly embedded in $ H=L^{2}(\Omega) $, hence, $ u_{m} $ converges strongly to $ u $ in $ L^{2}(0,T;H) $, which indicates $ u_{m}(0) $ converges $ u_{0} $.

\paragraph{}
Now, we introduce a lemma:
\begin{lemma}
Let $X$ be a given Banach space with dual $X'$ and let $u$ and $g$ be two functions belonging to $L^1(a,b;X)$. Then the following three conditions are equivalent:
\begin{enumerate}[label=(\roman*)]
    \item $u$ is almost everywhere equal to a primitive function of $g$, i.e., there exists $\xi\in X$ such that
    \[
    u(t)=\xi+\int_{a}^{t}g(s)ds, \quad \text{for a.e. } t\in[a,b].
    \]
    \item For every test function $\varphi\in\mathcal{D}(]a,b[)$,
    \[
    \int_{a}^{b}u(t)\varphi'(t)dt = -\int_{a}^{b}g(t)\varphi(t)dt \quad \left(\varphi'=\frac{d\varphi}{dt}\right).
    \]
    \item For each $\eta\in X'$,
    \[
    \frac{d}{dt}\langle u,\eta\rangle=\langle g,\eta\rangle
    \]
    in the scalar distribution sense on $]a,b[$.
\end{enumerate}
If $(i)-(iii)$ are satisfied we say that ${g\text{ is the }(X\text{-valued) distribution derivative of }u}$, and $u$ is almost everywhere equal to ${\text{a continuous function from }[a,b]\text{ into }X}$.\cite{temam2001}
\end{lemma}

Proof. Without loss of generality, let $[a, b]$ be $[0, T]$.

$(i)\Rightarrow(ii)$: Obtained by integration by parts; 
$(i)\Rightarrow(iii)$ is obvious.\par
$(iii)\Rightarrow(ii)$: For all $\varphi\in\mathcal{D}((0,T))$,
\[
\int_{0}^{T}\langle u(t),\eta\rangle\varphi'(t)dt = -\int_{0}^{T}\langle g(t),\eta\rangle\varphi(t)dt
\]
so we get, $\displaystyle \langle\int_{0}^{T}u(t)\varphi'(t)dt + \int_{0}^{T}g(t)\varphi(t)dt,\eta\rangle = 0$.\par
$(ii)\Rightarrow(i)$: Let $v = u - u_0$, where $u_0(t)=\int_{0}^{t}g(s)ds$. Obviously, $u_0(t)\in AC([0,T])$, so $u_0'(t) = g(t)$. From $(ii)$, for all $\varphi\in\mathcal{D}((0,T))$,
\[
\int_{0}^{T}v(t)\varphi'(t)dt = 0
\]
It is required to prove that $v(s)=\xi$ almost everywhere ($\xi \in X$).
Take $\varphi_0\in\mathcal{D}((0,T))$ such that $\int_{0}^{T}\varphi_0(t)dt = 1$. For all $\varphi\in\mathcal{D}((0,T))$, there exists a unique $\psi\in\mathcal{D}((0,T))$ such that $\varphi=\lambda\varphi_0+\psi'$, where $\lambda=\int_{0}^{T}\varphi(t)dt$.\par
Let $\xi=\int_{0}^{T}v(t)\varphi_0(t)dt\in X$. Next, prove that $v(t)=\xi$ almost everywhere.\par
It is easy to know that $\int_{0}^{T}(v(t)-\xi)\varphi(t)dt = 0$. Therefore, it only needs to prove that if $w\in L^1(a,b;X)$, for all $\varphi\in\mathcal{D}((0,T))$, if $\int_{0}^{T}w(t)\varphi(t)dt = 0$, then $w(t)=0$ almost everywhere for $t\in[0,T]$.\par

Define $\widetilde{w}(t)=
\begin{cases}
w(t),&t\in[0,T]\\
0,&\text{otherwise}
\end{cases}$, and take $\displaystyle \rho_{\varepsilon}=
\begin{cases}
C_{\varepsilon}\cdot e^{-\frac{1}{1 - |x|^2/\varepsilon^2}},&|x|\leq\varepsilon\\
0,&\text{otherwise}
\end{cases}$.Then
\[
\int_{0}^{T}w(t)(\rho_{\varepsilon}*\varphi)(t)dt=\int_{\mathbb{R}}\int_{\mathbb{R}}w(t)\rho_{\varepsilon}(z)\varphi(t - z)dzdt = 0
\]
\[
\int_{0}^{T}w(t)(\rho_{\varepsilon}*\varphi)(t)dt=\int_{\mathbb{R}}(\rho_{\varepsilon}*w)(t)\varphi(t)dt
\]
so $\int_{\mathbb{R}}(\rho_{\varepsilon}*w)(t)\varphi(t)dt = 0$.It is easy to verify that as $\varepsilon\rightarrow0$, $\rho_{\varepsilon}*\widetilde{w}\xrightarrow{L^1(\mathbb{R};X)}\widetilde{w}$.
Since $(\rho_{\varepsilon}*w)(t)$ is smooth, for all $\eta>0$, when $\varepsilon$ is sufficiently small, $(\rho_{\varepsilon}*w)(t)$ is $0$ on $[\eta,T - \eta]$.\\
Then we get $\widetilde{w}=0$ almost everywhere for $t\in[0,T]$.

\paragraph{}
Due to the weakly convergence of $ \frac{du_{m}}{dt}$ and lemma 2.6:
\begin{flalign}
    &\forall \varphi \in \mathscr{D}([0,T)) \\
    &(u,\varphi^{\prime}) = -(\frac{du_{m}}{dt},\varphi) \rightarrow -(\frac{du}{dt},\varphi) = (u,\varphi^{\prime}) \\
    &hence,\ u_{m}(t) \stackrel{V^{*}}{\longrightarrow} u(t) 
\end{flalign}
Furthermore,
\begin{equation}
    \left\{
    \begin{aligned}
        \lim_{m \to \infty}(\frac{d}{dt}u_{m},w_{j}) + v(Au_{m},w_{j}) + (B(u_{m}),w_{j}) &= \lim_{m \to \infty}(f,w_{j})\\
        \lim_{m \to \infty}u_{m}(0) &= \lim_{m\to \infty}P_{m}u_{0}
    \end{aligned}
    \right.
\end{equation}
\begin{equation}
    \Rightarrow\left\{
    \begin{aligned}
        \frac{d}{dt}u + vAu + B(u) &= f\\
        u(0) &= u_{0}
    \end{aligned}
    \right.
\end{equation}

\paragraph{}
Now, we need to regularize the solution space, illustrating $ u \in C([0,T];V) \cap L^{2}(0,T;D(A)) $. The only thing is to proof $ u_{m} \in L^{\infty}(0,T;V) \cap L^{2}(0,T;D(A)) $, because the following thing is similar to the above. We already have:
\begin{flalign}
\frac{d}{dt} \Vert u_{m} \Vert_{L^{2}(\Omega)}^{2} + v\lambda_{1}\Vert u_{m} \Vert_{V}^{2} \leq \frac{1}{v\lambda_{1}} \Vert f \Vert_{L^{2}(\Omega)}^{2}
\end{flalign}
\begin{flalign}
    \Rightarrow \int_{t}^{t+r} \Vert u_{m} \Vert_{V}^{2} ds \leq \frac{r}{v^{2}\lambda_{1}} \Vert f \Vert_{L^{2}}^{2} + \Vert u_{m} \Vert_{L^{2}}^{2}, \quad \Vert u_{m} \Vert_{L^{2}}^{2} \leq \Vert u_{m0} \Vert^{2} e^{-v\lambda_{1}t} + \frac{1}{v^{2}\lambda_{1}^{2}} \Vert f \Vert_{L^{2}}^{2} (1-e^{-v\lambda_{1}t})
\end{flalign}
Apply $ A_{m} $ to both sides of the equation:
\begin{flalign}
    (\frac{d}{dt}u_{m},Au_{m}) + v\Vert Au_{m} \Vert_{L^{2}} + (B(u_{m}),Au_{m}) &= (f,Au_{m})
\end{flalign}
\begin{flalign}
    \vert (B(u_{m}),Au_{m}) \vert &\leq \Vert u_{m} \Vert_{L^{\infty}} \Vert u_{m} \Vert_{V} \Vert Au_{m} \Vert_{L^{2}}\\
    &\leq c_{1} \Vert u_{m} \Vert_{L^{2}}^{1/2} \Vert u_{m} \Vert_{V}^{3/2} \Vert Au_{m} \Vert_{L^{2}}\\
    &\leq \frac{v}{4} \Vert Au_{m} \Vert_{L^{2}}^{2} + \frac{c_{1}^{\prime}}{v^{3}} \Vert u_{m} \Vert_{L^{2}}^{2} \Vert u_{m} \Vert_{V}^{4}
\end{flalign}
\begin{flalign}
    \Rightarrow \frac{d}{dt} \Vert u_{m} \Vert_{V}^{2} + v\Vert Au_{m} \Vert_{L^{2}}^{2} &\leq \frac{2}{v}\Vert f \Vert_{L^{2}} + \frac{2c_{1}^{\prime}}{v^{3}}\Vert u_{m} \Vert_{L^{2}} \Vert u_{m} \Vert_{V}^{4}\\
    and\ we\ have\ \Vert u_{m} \Vert_{V} &\leq \lambda_{1}^{-\frac{1}{2}} \Vert Au_{m} \Vert_{L^{2}}
\end{flalign}
\begin{lemma}[The Uniform Gronwall Lemma]\cite{temam2015}
    Let $ g $,$ h $,$ y $, be three positive locally integrable functions on $ [t_{0},+\infty] $ such that $ y^{\prime} $ is locally integrable on $ [t_{0},+\infty] $, and which satisfy:
    $$
    \frac{dy}{dt} \leq gy + h \quad \forall t\geq t_{0}
    $$
    $$
    \int_{t}^{t+r} g(s)ds \leq a_{1},  \int_{t}^{t+r} h(s)ds \leq a_{2},  \int_{t}^{t+r} y(s)ds \leq a_{3} \quad \forall t \geq t_{0},
    $$
    where $r,a_{1},a_{2},a_{3}$ are positive constants. Then
    $$
    y(t+r) \leq (\frac{a_{3}}{r}+a_{2})e^{a_{1}}, \quad \forall t \geq t_{0}.
    $$
\end{lemma}
Let $g,h,y$ respectively equals $ \frac{2c_{1}^{\prime}}{v^{3}}\Vert u_{m} \Vert_{L^{2}}^{2} \Vert u_{m} \Vert_{V}^{2}, \frac{2}{v}\Vert f \Vert_{L^{2}}^{2}, \Vert u_{m} \Vert_{V}^{2} $, it's obvious that $ a_{1},a_{2},a_{3} $ are bounded. Hence,
$$
\Vert u_{m} \Vert_{V}^{2} \leq (\frac{a_{3}}{r}+a_{2})e^{a_{1}},\quad when \ t\geq r(\forall r>0)
$$
Because $ u_{m0} \in D(A) $, $ u_{m} \in L^{\infty}(0,T;V) \cap L^{2}(0,T;D(A)). $ 

\paragraph{}
The last thing is to proof the uniqueness of the solution, the lemma we need as follows:
\begin{lemma}
Let $ V $, $ H $,$ V^{*} $ be three Hilbert spaces, each space included and dense in the following one, $ V $ being the dual of $ V $. If a function u belongs to $ L^2(0,T;V) $ and its derivative $ u' $ belongs to $ L^{2}(0,T;V) $, then $ u $ is almost everywhere equal to a function continuous from $ [0,T] $ into $ H $ and we have the following equality which holds in the scalar distribution sense on $ (0,T) $:
$$
\frac{d}{dt} \vert u \vert^{2} = 2\langle u',u \rangle.
$$\cite{temam2001}
\end{lemma}
Assuming $ u_{1},u_{2} $ are solutions of the equation, let $ u = u_{1}-u_{2} $:
\begin{equation}
    \left\{
    \begin{aligned}
        \frac{d}{dt}u + vAu &= B(u_{2})-B(u_{1})\\
        u(0) &= 0
    \end{aligned}
    \right.
    ,\quad u \in L^{\infty}(0,T;V) \cap L^{2}(0,T,D(A)) 
\end{equation}
By lemma 2.8:
$$
\frac{1}{2}\frac{d}{dt}\Vert u_{m} \Vert_{L^{2}}^{2} + va(u,u) = (B(u_{2})-B(u_{1}),u_{2}-u_{1})
$$
Then we observe:
$$
B(u_{1})-B(u_{2}) = B(u_{1},u)+B(u,u_{2})
$$
and
\begin{equation}
    \left\{
    \begin{aligned}
        (B(u_{1},u),u) &= 0\\
        \vert (B(u,u_{2}),u) \vert &= \vert \int_{\Omega} u \cdot \nabla u_{2} \cdot u dxdy\vert \leq \Vert u \Vert_{L^{\infty}}\Vert \nabla u_{2} \Vert_{V} \Vert u \Vert_{L^{2}}
    \end{aligned}
    \right.
\end{equation}
Due to $ H^{1}(\Omega) \hookrightarrow L^{\infty}(\Omega) $, with Agmon's inequality, $ \vert (B(u,u_{2}),u) \vert \leq c \Vert u \Vert_{L^{2}}^{3/2}\Vert u_{2} \Vert_{V}\Vert u \Vert_{V}^{\frac{1}{2}} $. By Young inequality, $ \vert (B(u,u_{2}),u) \vert \leq \frac{v}{2}\Vert u \Vert_{V}^{2} + \frac{c^{2}}{2v}\Vert u \Vert_{H}^{2} \Vert u_{2} \Vert_{V}^{2} $. Hence,
$$
\frac{d}{dt}\Vert u \Vert_{H}^{2} \leq \frac{c^{2}}{v}\Vert u \Vert_{H}^{2}\Vert u_{2} \Vert_{V}^{2}
$$
By Gronwall inequality, we obtain:
$$
\Vert u(t) \Vert_{H}^{2} \leq \Vert u(0) \Vert_{H}^{2} e^{\frac{c^{2}}{v}\int_{0}^{t}\Vert u_{2}(s) \Vert_{V}^{2}ds} = 0 
$$
Hence, the solution is unique. And the continuous dependence of the solution on $ f $ and the initial value $ u_0 $ can be deduced by the prior estimation.

\section{Application}
\paragraph{}
The Burgers equation, originally derived as a simplified model of fluid dynamics, has found widespread applications in various fields, including statistical mechanics, turbulence modeling, and transportation systems. Its nonlinearity and ability to describe shock formation and energy dissipation make it a powerful tool for analyzing complex systems. In transportation, the Burgers equation is particularly useful for modeling traffic flow dynamics, where it helps capture phenomena such as congestion, wave propagation, and the formation of traffic jams.

\subsection{Introduction of Detailed Application in Transportation Systems}
\paragraph{}
The application of the Burgers equation in transportation systems has been extensively studied. For instance, Payne (1971, 1979) utilized the Burgers equation to develop mathematical models for traffic flow, focusing on the propagation of shock waves and the behavior of traffic under varying conditions. These models provided insights into the dynamics of traffic congestion and helped improve traffic management strategies.\par

Kühne (1984) further expanded on this work by applying the Burgers equation to simulate traffic flow in urban networks. His studies demonstrated how the equation could be used to predict the formation and dissipation of traffic jams, offering valuable tools for urban planning and traffic control.\par

The Burgers equation's ability to describe the nonlinear interactions between vehicles and the environment makes it a cornerstone in traffic flow theory. By incorporating viscosity and external forcing terms, researchers can model the effects of road conditions, driver behavior, and traffic regulations on flow dynamics.

\subsection{Specific Process}
\paragraph{}
The standard form of the Burgers equation in traffic flow is:
\[
\frac{\partial u}{\partial t}+u\frac{\partial u}{\partial x}=\nu\frac{\partial^{2}u}{\partial x^{2}} + f(x,t)
\]
where:\par
$\displaystyle u(x,t)$ usually represents the speed or density of the traffic flow. If the density $\rho$ is taken as the variable, the equation needs to be transformed through the speed-density relationship, such as $\rho u = Q$ based on the conservation of flow. $\frac{\partial u}{\partial t}$ represents the rate of change of the traffic state with respect to time, reflecting the dynamic evolution of the traffic flow.

$\displaystyle u\frac{\partial u}{\partial x}$, the non-linear convection term, describes the "self-acceleration" effect of the vehicle density or speed. For example, high-density areas, that is, congestion points, will propagate upstream, similar to shock waves in fluids.

$\displaystyle \nu\frac{\partial^{2}u}{\partial x^{2}}$, the viscous diffusion term, simulates the smooth response of drivers to the following behavior. The larger $\nu$ is, the more gently the vehicles adjust their speeds, suppressing the formation of shock waves. $\nu\propto$ the reciprocal of the driver's reaction time. Conservative driving, that is, when $\nu$ is relatively large, will delay the formation of congestion. If $\nu\rightarrow0$, the equation degenerates into the inviscid Burgers equation, and shock waves will appear in the solution, corresponding to aggressive driving scenarios such as frequent lane-changing and emergency braking. If $\nu$ increases, for example, by introducing an adaptive $\nu(u)$, the lag effect of congestion propagation can be simulated, that is, the traffic recovery is slower than the deterioration.

$\displaystyle f(x,t)$, the external source term, represents external disturbances, such as traffic lights (periodic forcing terms), accidents (local pulse terms), or ramp traffic (spatially-dependent terms). For example, $f(x,t)=-\alpha\rho+\beta$ can represent ramp traffic, where $\beta$ is the inflow rate, or speed limit control ($\alpha$ is the speed adjustment coefficient). And adding the Dirac function $\delta(x - x_{k})$ can precisely describe the periodic influence of traffic lights at the intersection $x_{k}$.

In the previous sections, we have proven the existence of solutions under certain conditions. In traffic, this corresponds to the predictability of the traffic state, that is, given the initial conditions and boundary conditions, such as the inlet flow and the outlet capacity, the model can predict the traffic evolution. The uniqueness of the solution avoids prediction ambiguities. If $\nu = 0$, multi-valued solutions may occur at shock waves, corresponding to the observed phenomenon in traffic that "congestion waves can only propagate upstream". Below is the correspondence between the properties of the solutions of the Burgers equation and the actual traffic scenarios.

\begin{itemize}
    \item Shock solution: The solution has a discontinuity. For example, the density $\rho(x,t)$ suddenly changes at a certain position $x_{s}(t)$, which can explain the upstream propagation of the congestion front caused by traffic accidents. For example, the density behind the accident point $x_{s}(t)$ on the highway increases sharply, forming a "moving bottleneck".
    \item Rarefaction wave: The solution expands continuously, and the density gradient gradually decreases. It is used to simulate the dispersion of the traffic flow when the green light turns on or the restoration of the traffic flow downstream of the congestion.
    \item Viscous smoothing effect: The viscous term $\nu\frac{\partial^{2}u}{\partial x^{2}}$ will suppress the steepness of the shock waves, making the density change more gradual. In practice, this corresponds to the behavior of drivers "slowing down in advance and maintaining a safe distance", avoiding "phantom traffic jams", that is, congestion without an obvious cause.
\end{itemize}

The Burgers equation captures the spontaneous instability characteristics of the traffic flow through the non-linear convection term, the viscous term characterizes the smoothing effect of driving behavior on congestion, and the external term integrates control measures. The existence of its solutions ensures the reliability of the model. Shock waves and rarefaction waves correspond to the core dynamics of congestion formation and dissipation, and parameter adjustment can quantify the effects of management strategies. This provides a unified mathematical model basis for both micro-scale following behavior and macro-scale road network optimization. \cite{bec2007}

\bibliographystyle{plain}

\bibliography{mybib}

\end{document}